\documentclass[leqno,11pt]{amsart}
\usepackage{amsmath}
\usepackage[latin1]{inputenc}
\usepackage{amssymb}
\usepackage{amsmath}
\usepackage{amsfonts}
\usepackage{latexsym}
\usepackage[T1]{fontenc}
\usepackage[latin1]{inputenc}
\usepackage{epsfig}
\usepackage{graphics}
\usepackage{hyperref}
\usepackage[all,ps]{xy}
\RequirePackage{amsthm} \RequirePackage{amssymb}
\usepackage{latexsym}
\usepackage{amscd}
\usepackage{color}
\newtheorem{definition}{Definition}[section]
\newtheorem{lemma}[definition]{Lemma}
\newtheorem{prop}[definition]{Proposition}
\newtheorem{thm}[definition]{Theorem}
\newtheorem{cor}[definition]{Corollary}
\newtheorem{rmk}[definition]{Remark}

\def\dim{\mathop{\rm dim\, }\nolimits}
\def\proof{
  \noindent
  {\bf Proof: }}
\def\proofth{
  \noindent
  {\bf Proof of Theorem }}

\def\endproof{
{\unskip\nobreak\hfil\penalty50\hskip2em\hbox{}\nobreak\hfill
          $\square$\bigbreak}}

\begin{document}

\title{Stochastic Characterization of Harmonic maps on Riemannian polyhedra.}
\author{M. A. Aprodu, T. Bouziane}
\date{\today}
\maketitle

\begin{abstract}
The aim of this paper is to relate the theory of harmonicity, 
in the sense of Korevaar-Schoen and Eells-Fuglede, to the notion 
of Brownian motion in a Riemannian polyhedron, achieved by the 
second author. We define an exponential map at some 
singular points. Under the assumption that these
exponential maps are totally geodesic (for instance
in dimension one), we find the infinitesimal generator 
of a Brownian motion in a Riemannian polyhedron. We prove 
that it is uniquely defined on some Banach
space. Finally, we show that harmonic maps, in the sense of 
Eells-Fuglede, with target smooth Riemannian manifolds, are 
characterized by mapping Brownian motions in Riemannian
polyhedra into martingales, while harmonic morphisms are exactly the 
maps which are Brownian preserving paths.
\end{abstract}

\section{Introduction.}

Brownian motions in Riemannian manifolds are intimately related to harmonic functions, maps
and morphisms. The origin of this relationship is the definition of a Brownian motion in a
Riemannian manifold as a diffusion process generated by the Laplace-Beltrami operator, which
is also the basic tool in the theory of harmonic maps. In \cite{D}, Darling studied the
relation between the behaviour of Brownian motions under maps between Riemannian manifolds
and harmonicity.

The theory of harmonic maps between smooth Riemannian manifolds was extended by Gromov,
Korevaar and Schoen (see \cite{GS}, \cite{KS}) to the case of maps between certain singular
spaces, such as admissible Riemannian polyhedra. Riemannian polyhedra are interesting in
many regards. They carry natural harmonic space structures in the sense of Brelot, and they
include several geometric objects: smooth Riemannian manifolds, Riemannian orbit spaces,
normal analytic spaces, Thom spaces etc. A deeper study of harmonic maps and morphisms
between Riemannian polyhedra was done by Eells and Fuglede in \cite{EF}. Using the {\em weak
conformality property}, in the case when the target space is a Riemannian manifold, they
obtained the same characterization for harmonic morphisms as in the smooth case, cf.
\cite{F} and \cite{I}.

On the other hand, a rigorous construction of a Brownian motion in a Riemannian complex, was
given by the second author in \cite{TB}. In particular, this construction applies to the
case of a Riemannian polyhedron. Previously, Brin-Kifer had constructed a Brownian motion in
the particular case of flat $2$-dimensional admissible complexes \cite{BK}.

The aim of this paper is to relate, in the case of Riemannian polyhedra, the theory of
harmonic maps and morphisms \cite{EF} to the notion of Brownian motions in Riemannian
polyhedra \cite{TB}, in order to generalize Darling's results, \cite{D}, \cite{M}. Notice
that the second differential calculus on Riemannian manifolds, which is the basis of the
theory of stochastic calculus, has no natural generalization on Riemannian polyhedra.
Consequently, we are compelled to develop a new approach mixing smooth theory with some
hybrid methods.

The outline of the paper is as follows.  Section 2, included here for the sake of
completeness, is an overview on Riemannian polyhedra, energy of maps, harmonic maps and
morphisms on Riemannian polyhedra, Brownian motions in  Riemannian manifolds, martingales
etc. In Section 3, we prove that the Brownian motion has a unique infinitesimal generator
defined on some Banach space. In this section we also study the behaviour of Brownian
motions under harmonic functions in the sense of Gromov-Korevaar-Schoen \cite{GS}, \cite
{KS}. In order to state this characterization we show that the Brownian motion in a
Riemannian polyhedron has the "Laplacian" as an infinitesimal generator Theorem \ref{gen}
(we give a suitable definition of the "Laplacian" in this case). 
For technical reasons we make some assumption 
on the exponential maps (being totally geodesic). This condition is
realized for example in the one-dimensional case
or for flat metrics.
In the last section we prove that harmonic maps, with target smooth Riemannian manifolds in the sense of \cite{EF},
are exactly those which map Brownian motions in Riemannian polyhedra into martingales,
Theorem \ref{th2}, while harmonic morphisms are exactly the maps which are Brownian
preserving paths, Theorem \ref{th3}.

\medskip

\noindent{\bf Acknowledgements.} The authors would like to thank Professor A. Verjovsky to
encourage them to investigate these problems and the ICTP Trieste for hospitality during
this work. The first named author was partially supported by the ANCS contract
2-CEx06-11-12/25.07.2006.

\section{Preliminaries.}

In this section we recall some basic notions and results which will be used throughout the
paper.

\subsection{Riemannian admissible complexes} \cite{BB}, \cite{Br},
\cite{BH}, \cite{DJ}, \cite{EF}, \cite{T}.

Let $C$ be a locally finite simplicial complex, endowed with a piecewise smooth Riemannian
metric $g$,  i.e. $g$ is a family of smooth Riemannian metrics $g_S$ on simplexes $S$ of
$C$, such that the restriction ${(g_S)}_{|S'}=g_{S'}$, for any simplexes $S'$ and $S$ with
$S'\subset S$.

The set of all formal linear combinations $\alpha=\sum_{v\in C} \alpha(v)v$ of vertices of
$C$, such that $0\leqslant \alpha(v)\leqslant 1$, $\sum_{v\in C} \alpha(v)=1$ and $\{v;
\alpha(v)>0\}$ is a simplex of $C$, is denoted by $space |C|$. This set is a subset of the
linear space  $lin C$ of all formal finite linear combinations of vertices of $C$. A vertex
$v$ of $C$ will be identified with the formal linear combination $1v$, thus formal linear
combinations of vertices become true linear combinations in $lin C.$

Let $C$ be a finite dimensional simplicial complex which is connected locally finite. A map
$f$ from $[a,b]$ to $C$ is called a broken geodesic if there is a subdivision
$a=t_0<t_1<\dots<t_{p+1}=b$, such that $f([t_i,t_{i+1}])$ is contained in some cell and the
restriction of $f$ to $[t_i,t_{i+1}]$ is a geodesic inside that cell. Then, define the
length of the broken geodesic map $f$ to be:
$$
\mathcal{L}(f)=\sum_{i=0}^p d(f(t_i),f(t_{i+1})).
$$
The length inside each cell is measured with respect to its metric.

For every two points $x,y$ in $C$, define $\tilde d(x,y)$ to be the lower bound of the
lengths of broken geodesics from $x$ to $y$. Note that $\tilde d$ is a pseudo-distance.

If $C$ is connected and locally finite, then $(C,\tilde d)$ is a length space and hence a
geodesic space (i.e. a metric space where every two points are connected by a curve with
length equal to the distance between them ), if complete.

We say that the complex $C$ is {\em admissible}, if it is dimensionally homogeneous, and for
every connected open subset $U$ of $C$, the open set $U\setminus \{ U\cap \{(n-2)-\mbox{
skeleton }\} \}$ is connected, where $n$ is the dimension of $C$ (i.e. $C$ is
$(n-1)$-chainable).

We call an admissible  connected locally finite simplicial complex, endowed with a piecewise
smooth Riemannian metric, an {\em admissible Riemannian complex}.

In the sequel we shall denote a {\em p-skeleton} of a complex $C$ by $C^{(p)}$.

\subsection{Riemannian polyhedra} \cite{EF}.

We mean by {\em polyhedron} a connected locally compact separable Hausdorff space $K$ for
which there exists a simplicial complex $C$ and a homeomorphism $\theta : C \rightarrow K$.
Any such pair $(C,\theta )$ is called a {\em triangulation} of $K$. The complex $C$ is
necessarily countable and locally finite (see \cite{S} page 120) and the space $K$ is path
connected and locally contractible. The {\em dimension} of $K$ is by definition the
dimension of $C$ and it is independent of the triangulation.

If $K$ is a polyhedron with specified triangulation $(C,\theta)$, we shall speak of
vertices, simplexes, $i$-skeletons (the set of simplexes of dimensions lower or equal to
$i$), stars of $K$ as the image under $\theta$ of vertices, simplexes, $i$-skeletons, stars
of $C$. Thus our simplexes become compact subsets of $K$.

If for a given triangulation $(C,\theta)$ of the polyhedron $K$, the homeomorphism $\theta$
is locally bi-lipschitz then $K$ is said to be a {\em Lip polyhedron} and $\theta$ a {\em
Lip homeomorphism}.

A {\em null set} in a Lip polyhedron $K$ is a set $Z\subset K$ such that $Z$ meets every
maximal simplex $S$, relative to a triangulation $(C,\theta)$ (hence any) in a set whose
pre-image under $\theta$ has $n$-dimensional Lebesgue measure $0$, with $n=\dim S$. Note
that {\em 'almost everywhere' } (a.e.) means everywhere except in some null set.

A {\em Riemannian polyhedron} $K=(K,g)$  is defined as a Lip polyhedron $K$ with a specified
triangulation $(C,\theta)$, such that $C$ is a simplicial complex endowed with a covariant
bounded measurable Riemannian metric tensor $g$, satisfying the ellipticity condition below.
In fact, suppose that $K$ has homogeneous dimension $n$ and choose a measurable Riemannian
metric $g_S$ on the open Euclidean $n$-simplex $\theta^{-1}(\; \mathop{S}\limits^o\; )$ of
$C$. In terms of Euclidean coordinates $\{x_1,\dots,x_n\}$ of points $x=\theta^{-1}(p)$,
$g_S$ assigns to almost every point $p\in \mathop{S}\limits^o$ (or $x$), an $n\times n$
symmetric positive definite matrix $g_S = (g_{ij}^S(x))_{i,j=1,\dots,n}$, with measurable
real entries and there is a constant $\Lambda_S >0$ such that (ellipticity condition):
$$
\Lambda_S^{-2}\sum_{i=0}^n(\xi^i)^2\leq \sum_{i,j}
g^S_{ij}(x)\xi^i\xi^j\leq\Lambda_S^2\sum_{i=0}^n(\xi^i)^2
$$
for $a.e.$ $x\in\theta^{-1}(\; \mathop{S}\limits^o\; )$ and every $\xi=(\xi^1,\dots,\xi^n)
\in \mathbb{R}^n$. This condition amounts to the components of $g_S$ being bounded and it is
independent not only of the choice of the Euclidean frame on $\theta^{-1}(\;
\mathop{S}\limits^o\; )$ but also of the chosen triangulation.

For simplicity of statements, we shall sometimes require that, relative to a fixed
triangulation $(C,\theta)$ of the Riemannian polyhedron $K$ (uniform ellipticity condition),
$$
\Lambda:= \mbox{ sup } \{\Lambda_S:S \mbox{ is simplex of } K\}<\infty.
$$

A Riemannian polyhedron $K$ is said to be {\em admissible} if for a fixed triangulation
$(C,\theta)$ (hence any) the Riemannian simplicial complex $C$ is admissible.

We underline that, for simplicity, the given definition of a Riemannian polyhedron $(K,g)$
contains already the fact (because of the definition above of the Riemannian admissible
complex) that the metric $g$ is {\em continuous} relative to some (hence any) triangulation
(i.e. for every maximal simplex $S$ the metric $g_S$ is continuous up to the boundary). This
fact is sometimes omitted in the literature. The polyhedron is said to be {\em simplexwise
smooth } if relative to some triangulation $(C,\theta)$ (and hence any), the complex $C$ is
simplexwise smooth. Both continuity and simplexwise smoothness are preserved under
subdivision.

\subsection{Energy of maps} \cite{GS}, \cite{KS}, \cite{EF}.

The concept of energy of maps from a Riemannian domain into an arbitrary metric space $Y$
was defined and investigated by Gromov, Korevaar and Schoen \cite{GS}, \cite{KS}. Later on,
this concept was extended by Eells and Fuglede \cite{EF} to the case of a map from an
admissible Riemannian polyhedron $K$ with simplexwise smooth Riemannian metric. The {\em
energy} $E(\varphi)$ of a map $\varphi$ from $K$ to the space $Y$ is defined as the limit of
suitable approximate energy expressed in terms of the distance function $d_Y$ of $Y$.

The maps $\varphi : K\rightarrow Y$ of finite energy are precisely those quasicontinuous
(i.e. have  continuous restrictions to closed sets), whose complements have arbitrarily
small capacity, (see \cite{EF}, page 153) whose restriction to each top dimensional simplex
of $K$ has finite energy in the sense of Korevaar-Schoen, and $E(\varphi)$ is the sum of the
energies of these restrictions.

Consider now an admissible $m$-dimensional Riemannian polyhedron $(K,g)$ with simplexwise
smooth Riemannian metric. It is not required that $g$ is continuous across lower dimensional
simplexes. The target $(Y,d_Y)$ is an arbitrary metric space.

Denote $L^2_{loc}(K,Y)$ the space of all $\mu_g$-measurable ($\mu_g$ the volume measure of
$g$) maps $\varphi :K\rightarrow Y$ having separable essential range and for which the map
$d_Y(\varphi (.),q)\in L^2_{loc}(K,\mu_g)$ (i.e. locally $\mu_g$-squared integrable) for
some point $q$ (hence by triangle inequality for any point). For $\varphi,\psi \in
L^2_{loc}(K,Y)$ define their distance $\mathcal{D}(\varphi,\psi)$ by:
$$
 \mathcal{D}^2(\varphi,\psi)= \int\limits_K d_Y^2(\varphi (x),\psi(y)) d\mu_g(x).
$$

Two maps $\varphi,\psi \in L^2_{loc}(K,Y)$ are said to be {\em equivalent} if
$\mathcal{D}(\varphi,\psi)=0$, (i.e. $\varphi(x)=\psi(x)$ $\mu_g$-a.e.). If the space $K$ is
compact, then $\mathcal{D}(\varphi,\psi)<\infty$ and $\mathcal{D}$ is a metric on
$L^2_{loc}(K,Y)=L^2(K,Y)$ which is complete if the space $Y$ is complete \cite{KS}.

The {\em approximate energy density} of the map $\varphi\in L^2_{loc}(K,Y)$ is defined for
$\epsilon >0$ by:
$$
 e_\epsilon(\varphi)(x)=\int\limits_{B_K(x,\epsilon)}\frac{d_Y^2(\varphi(x),\varphi(x'))}
 {\epsilon^{m+2}}d\mu_g(x').
$$
The function $e_\epsilon(\varphi)\ge 0$ is locally $\mu_g$-integrable.

The {\em energy} $E(\varphi)$ of a map $\varphi$ of class $L^2_{loc}(K,Y)$ is:
$$
 E(\varphi)=\sup_{f\in
 \mathcal{C}_c(K,[0,1])}\left(\limsup_{\epsilon\rightarrow 0}\int\limits_K
 fe_\epsilon(\varphi) d\mu_g\right),
$$
where $\mathcal{C}_c(K,[0,1])$ denotes the space of continuous functions from $K$ to the
interval $[0,1]$ with compact support.

A map $\varphi: K\rightarrow Y$ is said to be {\em locally of finite energy}, and we write
$\varphi \in W^{1,2}_{loc}(K,Y)$, if $E(\varphi_{|U})<\infty$ for every relatively compact
domain $U\subset K$, or, equivalently, if $K$ can be covered by domains $U\subset K$ such
that $E(\varphi_{|U})<\infty$.

For example (Lemma 4.4, \cite{EF}), every Lip continuous map $\varphi: K \rightarrow Y$ is
of class $W^{1,2}_{loc}(K,Y)$. In the case when $K$ is compact, $W^{1,2}_{loc}(K,Y)$ is
denoted by $W^{1,2}(K,Y)$ the space of all maps of finite energy.

$W^{1,2}_c(K,Y)$ denotes the linear subspace of $W^{1,2}(K,Y)$ consisting of all maps of
finite energy of compact support in $K$.

\subsection{Harmonic maps and harmonic morphisms on Riemannian polyhedra} \cite{EF}.

Let $(K,g)$ be an arbitrary admissible Riemannian polyhedron ($g$ is only bounded,
measurable, with local elliptic bounds), $\dim K=m$ and $(Y,d_Y)$ a metric space .

A continuous map $\varphi:K\rightarrow Y$ of class $W_{loc}^{1,2}(K,Y)$ is said to be {\em
harmonic} if it is {\em bi-locally E-minimizing}, i.e. $K$ can be covered by relatively
compact subdomains $U$ for each of which there is an open set $V\supset \varphi(U)$ in $Y$
such that
$$
E(\varphi_{|U})\leq E(\psi_{|U})
$$
for every continuous map $\psi\in W_{loc}^{1,2}(K,Y)$, with $\psi(U)\subset V$ and
$\psi=\varphi$ in $K\backslash U$.

Let $(N,h)$ denote a smooth Riemannian manifold without boundary of dimension $n$ and
$\Gamma^k_{\alpha\beta}$ the Christoffel symbols on $N$. By a {\em weakly harmonic map }
$\varphi:K\rightarrow N$ we mean  a quasicontinuous map (a map which is continuous on the
complement of open sets of arbitrarily small capacity; in the case of the Riemannian
polyhedron $K$, it is just the complement of open subsets of the $(m-2)$-skeleton of $K$) of
class $W_{loc}^{1,2}(K,N)$ with the following property:

For any chart $\eta:V\rightarrow \mathbb{R}^n$ on N and any quasiopen set
$U\subset\varphi^{-1}(V)$ of compact closure in $K$, the equation
$$
\int\limits_U\langle\nabla\lambda,\nabla\varphi^k\rangle d\mu_g= \int\limits_U
\lambda(\Gamma^k_{\alpha\beta}\circ\varphi)
\langle\nabla\varphi^\alpha,\nabla\varphi^\beta\rangle d\mu_g,
$$
holds for every $k=1,\dots,n$ and every bounded function $\lambda\in W_0^{1,2}(U)$.

When $K$ and $Y$ denote two Riemannian polyhedra (or any harmonic spaces in the sense of
Brelot; see Chapter 2, \cite{EF}), a continuous map $\varphi:K\rightarrow Y$ is a {\em
harmonic morphism} if, for every open set $V\subset Y$ and for every harmonic function $v$
on $V$, $v\circ\varphi$ is harmonic on $\varphi^{-1}(V)$.

\subsection{Brownian motions in Riemannian manifolds} \cite{D},  \cite{ME}, \cite{V}.

Consider $(\Omega, \mathcal{A}, P)$ a probability space, $(E, \varepsilon)$ a measurable
space, and $I$ an ordered set. By a {\em stochastic process} on $(\Omega, \mathcal{A}, P)$
with values on $(E, \varepsilon)$ and $I$ as time interval, we mean a map (see \cite {ME},
or \cite{V}, or \cite{D}):
$$
\begin{array}{cccc}
  X: & I\times\Omega & \rightarrow & E \\
   & (t,\omega)& \mapsto & X( t,\omega),\\
\end{array}
$$
such that for each  $t\in I$, $X_t: \omega\in \Omega \mapsto X(t,\omega)\in  E$ is
measurable from $(\Omega,\mathcal{A})$ to $(E,\varepsilon)$.

A family $\mathcal{F}=(\mathcal{F}_t)_{t\in I}$ of $\sigma$-subalgebras of $\mathcal{A}$,
such that $\mathcal{F}_s\subset\mathcal{F}_t$, for all $s$, $t$ with $s< t$, is called {\em
a filtration } on $(\Omega, \mathcal{A}, P)$ with $I$ time interval.

Given a filtration $\mathcal{F}=(\mathcal{F}_t)_{t\in I}$, a process $X$, admitting as time
interval a part $J$ of $I$, is said to be {\em adapted } to $\mathcal{F}$, if for every
$t\in J$, $X_t$ is $\mathcal{F}_t$-measurable.

A real-valued process $X$ is said to be  a {\em submartingale}, with respect to a filtration
$\mathcal{F}_t$ fixed on $(\Omega, \mathcal{A}, P)$, if it has the following properties : a)
$X$ is adapted; b) each random variable $X_t$ is integrable; c) for each pair of real
numbers $s$, $t$, $s<t$, and every $A\in \mathcal{F}_s$ we have:
$$
\int\limits_A X_sdP\leq\int\limits_AX_tdP.
$$
When the equality holds we say that $X$ is a {\em martingale}.

A real-valued process $X$ is said to be a {\em continuous local martingale} if and only if
it is a continuous ( with respect to the time variable ) adapted process $X$ such that each
$X_{t\wedge T_n}\chi_{\{T_n>0\}}$ is a martingale, where $\chi$ is the characteristic
function and $T_n$ is the stopping time: $inf\{t:|X_t|\geq n\}$.

A {\em semimartingale} is the sum of a continuous local martingale and a process with finite
variation. If the process of the finite variation is an increasing one, the semimartingale
is called a {\em local submartingale}.

Let $M$ be a manifold with a connection $^M\nabla$, and $X$ be an $M$-valued process.
Following Schwartz characterization \cite{Sc}, a {\em $^M\nabla$-martingale tester},
$(U_1,U_2,U_3,f)$ will consist of:
\begin{itemize}
    \item open sets $U_1$, $U_2$, $U_3$ in $M$ with
    $\overline{U_1}\subset U_2 \subset \overline{U_2}\subset U_3,$
     \item a convex function $f:U_3\rightarrow \mathbb{R}.$
\end{itemize}
The process $X=(X_t,\mathcal{F}_t)$ is said to be a {\em $^M\nabla$-martingale}, if it is a
continuous semimartingale on $M$ (i.e. $\forall f\in \mathcal{C}^2(M)$, $f\circ X$ is a real
valued semimartingale), and for all $^M\nabla$-martingale tester $(U_1,U_2,U_3,f)$, the
process
$$Y=(Y_t,\mathcal{F}_t), \,
Y_t=\int\limits_0^t\chi_F(s)d(f\circ X_s),
$$
is a local submartingale. $F$ denote the previsible set $\bigcup\limits_{i=1}^\infty
(\sigma_i,\tau_i]$ where $\sigma_i,\tau_i, i\geq 0$ is the collection of stopping-times,
associated to the process $X$ and any $^M\nabla$-martingale tester, defined by:
\begin{itemize}
    \item $\sigma_0=0, \tau_0=0$
    \item $\sigma_i=inf\{t> \tau_{i-1}:X_t\in U_1\};
    \tau_i=inf\{t> \sigma_i:X_t\notin
    \overline{U_2}\}, i\geq 1.$
\end{itemize}
$\chi_F$ denotes the characteristic function.

Suppose $M$ is a Riemannian manifold with Levi-Civita connection $^M\nabla$. A Brownian
motion is characterized as a diffusion $B=(B_t,\mathcal{F}_t)$ with generator
$\frac{1}{2}\Delta$; in other words, for all $f:M\rightarrow \mathbb{R}$, the process $C^f$,
where $C^f_t=f(B_t)-f(B_0)-\frac{1}{2}\int\limits_0^t\chi_F(s)\Delta f(B_s)ds$, is a local
martingale.

\section{Brownian motions in admissible Riemannian  polyhedra.}

In \cite{TB}, the second author proved the existence of a Brownian motion on a Riemannian
polyhedron. In this section, we find explicitly its infinitesimal generator, Theorem
\ref{gen}. Furthermore, we give necessary and sufficient conditions for harmonicity in terms
of local martingales, see Corollary \ref{th1}.

\medskip

We begin by recalling some basic results used in the sequel.

\subsection{Boundary normal coordinates, \cite{KKL}, \cite{SU}
and the second fundamental form} \cite{BW}, \cite{ME}, \cite{NS}. 
\label{sec}

Let $M$ be a Riemannian manifold with non-empty boundary 
$\partial M$. For any point $y\in M$ there is a shortest 
geodesic to the boundary that is normal to $\partial M$.

Similarly to the exponential mapping defined on $T_xM$ for 
$x$ an interior point of $M$, we can define (see \cite{KKL}) 
the {\em boundary exponential mapping}: $\mbox{exp}_{\partial
M}:\partial M\times \mathbb{R}_+ \rightarrow M$, 
$\mbox{exp}_{\partial M}(z,t)=\gamma_{z,\nu}(t)$, 
where $\mathbb{R}_{+}=[0,\infty)$ and $t$ sufficiently small
such that $\mbox{exp}_{\partial M}(z,t)\in M.$ Here, 
$\gamma_{z,\nu}$ denotes the normal geodesic to $\partial M$ 
whose derivative at zero equals $\nu$, the unitary normal vector to
$\partial M$ at the point $z$.

Using the boundary exponential mapping, one introduced (see for example \cite{KKL} or
\cite{SU}) the {\em boundary normal} (or {\em semi-geodesic}) {\em coordinates}, analogously
to the Riemann normal coordinates. Compared to the classical case of empty boundary, instead
of a set of geodesics starting  from a point one considers the set of geodesics normal to
$\partial M.$

Consider $\mathcal{U}_{\rho}=\partial M \times [0,\rho)$ a collar neighbourhood of
$\partial M \times \{0\}$ in the boundary cylinder $\partial M \times \mathbb{R}_+$. Denote
by
$$
\mathcal{V}_{\rho}=\mbox{exp}_{\partial M}(\mathcal{U}_{\rho})= \{x\in M; d(x, \partial M)<
\rho\}
$$
a collar neighbourhood of $\partial M$ in $M$.

Then, for $\rho$ sufficiently small, define $(\mathcal{V}_{\rho}, x_1,\dots,x_n)$ local
coordinates  in $M$ ( the boundary normal coordinates) in the following way: for $x\in
\mathcal{V}_{\rho}$, $x_n:=d(x, \partial M)$, $z\in \partial M$ is the unique boundary point
such that $d(x,z)=d(x, \partial M)$ and $(x_1,\dots,x_{n-1})$ on $\partial M$ are local
coordinates around $z$. $\rho$ is chosen small enough such that $\gamma_{z,\nu}(t)$ is  the
unique shortest geodesic to $\partial M$ for $t< \rho$.

As in the case of Riemannian manifolds without boundary, the Laplacian  is given in boundary
normal coordinates (see \cite{L}) by:
$$
\Delta = \frac{1}{\sqrt{|g|}}\frac{\partial}{\partial x_i}
              \left( g^{ij}\sqrt{|g|}\frac{\partial}{\partial x_j}\right)=
              g^{ij}\left( \frac{{\partial}^2}{\partial x_i\partial x_j}-
              \Gamma^k_{ij}\frac{\partial}{\partial x_k}\right) ,
$$
where $|g|=det(g_{kl})$ and $\Gamma^k_{ij}$ are the Christoffel symbols on $M$.

\medskip

Let $(M,g)$, $(N,h)$ two Riemannian manifolds and $\varphi: M\rightarrow N$ a smooth map.
The Levi-Civita connection $\nabla^M$ of $M$ and the pull-back connection $\nabla^{\varphi}$
of the pull-back bundle $\varphi^{-1}TN$ induce a connection $\nabla$ on the bundle
$T^*M\otimes \varphi^{-1}TN$. Applying this connection to $d\varphi$ one obtain the {\em
second fundamental form of $\varphi$} (also called the {\em Hessian of $\varphi$}) \;(see
\cite{BW} or \cite{NS}), denoted by $\mbox{Hess }\varphi$, or $\nabla d\varphi$, and
explicitly given by:

\begin{equation}\label{hes1}
(\mbox{Hess }\varphi)(X,Y):=\nabla^{\varphi}_X(d\varphi(Y))- d\varphi(\nabla^M_X Y), \;
\forall X,Y\in \Gamma(TM).
\end{equation}

When $N=\mathbb{R}$ and $(M,g)$ is a manifold with or without boundary then (\ref{hes1})
reads:
\begin{equation}\label{hes2}
(\mbox{Hess }\varphi)(X,Y):=X(Y(\varphi))- d\varphi(\nabla^M_X Y), \; \forall X,Y\in
\Gamma(TM).
\end{equation}
Using the Hessian, the Laplacian can be also defined as:
$$
\Delta f:=Tr (\mbox{Hess }f).
$$

A smooth map $\varphi:M\rightarrow N$ between two Riemannian manifolds is {\em totally
geodesic} (see \cite{BW}, \cite{ME}) if for every $f:N\rightarrow \mathbb{R}$
$$
(\mbox{Hess}^M (f\circ \varphi))(*,*)=(\mbox{Hess}^N f)(d\varphi(*),d\varphi(*)),
$$
where $\mbox{Hess}^M$, $\mbox{Hess}^N$ denote the second fundamental forms on $M$, $N$
respectively.

\subsection{The tangent cone and the exponential map.}
\label{exp}

Let $(K,g)$ be an $n$-dimensional admissible Riemannian polyhedron 
and $p$ a point in the ($(n-1)$-skeleton)$\setminus$ ($(n-2)$-skeleton).

We shall  slightly reformulate the {\em definition of the tangent cone} 
previously introduced in \cite{BB}.

Suppose that $p$ is in ${\mathop{\widehat{S_{n-1}}}\limits^o},$ the topological interior of
the $(n-1)$-simplex $S_{n-1}$. Let $S_n^1,S_n^2,\dots,S_n^k$, $k\geq 2$, denote the
$n$-simplexes adjacent to $S_{n-1}$. Then each $S_n^l$, for $l=1,\dots,k$, can be viewed as
an affine simplex in $\mathbb{R}^n$, that is $S_n^l=\bigcap_{i=0}^n H_i$ where $H_i$ are
closed half spaces in $\mathbb{R}^n.$ The Riemannian metric $g_{S_n^l}$ is the restriction
to $S_n^l$ of a smooth Riemannian metric defined in an open neighbourhood of $S_n^l$ in
$\mathbb{R}^n$.

Since $p \in ((n-1)-skeleton)\setminus ((n-2)-skeleton)$, each 
$S_n^l$ for $l=1,\dots,k$, can be viewed, locally around $p$, as a 
manifold with boundary, where the boundary is
$S_{n-1}.$ Then  there exists a unique
hyperplane, for $i=0,\dots,n,$ containing $p$. Define 
$T_pS_n^l$ as the half-space $H_i$ which contains the 
corresponding hyperplane.

Notice that $T_pS_n^l$ can be naturally embedded in \;$lin S_n^l \subset lin K $ and
\begin{equation}\label{tan}
T_pS_n^l=T_pS_{n-1}\times [0,\infty).
\end{equation}

\medskip

Define the {\em tangent cone of $K$ over $p$} as: $T_pK=\bigcup_{l=1}^k T_pS_n^l\subset lin K.$

\medskip

The difference from the original definition (see \cite{BB}) is that we do not need to pass
to subdivision of $K$ in order to make the point $p$ become a vertex.

Let $\mathfrak{U}_p(S_n^l)$ be the subset of all unit vectors in $T_pS_n^l$ and denote
$\mathfrak{U}_p=\mathfrak{U}_p(K) =\bigcup_{S_n^l\ni p} \mathfrak{U}_p(S_n^l)$. The set
$\mathfrak{U}_p$ is called the {\em link of $p$ in $K$}. As $S_n^l$ is a simplex adjacent to
$p$, then $g_{S_n^l}(p)$ defines a Riemannian metric on the $(n-1)$-simplex
$\mathfrak{U}_p(S_n^l)$. The family $g_p$ of Riemannian metrics $g_{S_n^l}(p)$ turns
$\mathfrak{U}_p(S_n^l)$ into a simplicial complex with a piecewise smooth Riemannian metric
such that the simplexes are spherical.

\medskip

Having defined the tangent cone, and using the boundary exponential
map, we can introduce next the {\em exponential map} locally around
a point $p$ in the topological interior 
of an $(n-1)$-simplex $S_{n-1}$.

Take $V_0$ a neighbourhood of $0$ in
$T_pK.$ The definition of the exponential map $E_p : T_pK \rightarrow K$ on each maximal face
$V_0\bigcap T_pS_n^l$, $l=1,\dots,k$ is based on the fact that, locally around $p$, each
$S_n^l$ becomes a manifold with boundary, with $\partial(S_n^l)=S_{n-1}$. This allows us to
consider the boundary exponential map (see Section \ref{sec}):
$$
\mbox{exp}_{\partial  S_n^l }:\mathcal{V}_{\rho}\rightarrow \mathcal{W}_{\rho},
$$
where $U_p$ is a small neighbourhood of $p$, $\mathcal{V}_{\rho}=(U_p\cap S_{n-1})\times
[0,\rho)$ is a collar neighbourhood  of $(U_p\cap S_{n-1})\times {0}$ in the boundary
cylinder $(U_p\cap S_{n-1})\times \mathbb{R}_+$ and
$$
\mathcal{W}_{\rho}=\mbox{exp}_{\partial S_n^l }(\mathcal{V}_{\rho}):= \{x\in (U_p\cap S_n^l)
; d(x,(U_p\cap S_{n-1}) )< \rho\}.
$$

Moreover, on the manifold $S_{n-1}$ we consider the usual exponential map at $p$:
$$
\mbox{exp}_p:T_pS_{n-1}\rightarrow S_{n-1}.
$$

Using the decomposition $T_pS_n^l=T_pS_{n-1}\times [0,\infty)$, we
define the exponential map
$$
E_p : V_0\cap T_pS_n^l\rightarrow  S_n^l
$$
in the following way. Consider $u$ a tangent vector in $V_0\cap T_pS_n^l$. We can decompose
$u=(v,w)$ where $v\in T_pS_{n-1}$ and $w$ is a normal vector to $\partial  S_n^l$. Then
$$
E_p(u)= \mbox{exp}_{\partial S_n^l}(\mbox{exp}_p(v),||w||).
$$

\subsection{Brownian motions.}

The Brownian motion in a piecewise smooth Riemannian complex, was obtained in \cite{TB}, as
a weak limit of isotropic processes. This construction holds obviously for the piecewise
smooth Riemannian polyhedron $K$.

Let us recall some essential facts about this construction. In \cite{TB}, the second author
defined a process: $Y^\eta= (\Omega,\mathcal{F}^0_t,Y_t^\eta,\theta_t,P)$, for $\eta\in
(0,1]$, in the following way:
$$
Y_t^{\eta}(\omega) =\left\{
\begin{array}{ccc}
  \Upsilon_{\eta Z_i(\omega)}(\frac{t}{\eta^2}-\tau_i(\omega))& \mbox { if } &
  \tau_i(\omega)\leq \frac{t}{\eta^2}\leq \tau_{i+1}(\omega)\\
 D & \mbox { if } & \xi(\omega)\leq \frac{t}{\eta^2},\\
\end{array}\right.
$$
where $\tau_i$ are the stopping times such that, for all $i\in {\bf N}$, the real random
variable $(\tau_i-\tau_{i+1})$ is exponentially distributed and $\tau_0=0$; $\Upsilon_\eta$
is the generalized geodesic flow (see \cite{BB}); $D$ is the one point compactification of
$K$ (because $K$ is semicompact) and $\xi$ is the life time of $Y_t^{\eta}$; $Z$ is a unit
tangent vector randomly chosen in the link of the point $\Upsilon_{\eta
Z_{i-1}(\omega)}(\tau_i(\omega))$ with respect to the volume measure (link is viewed as a
spherical Riemannian polyhedron), where $Z_0(\omega)$ is also a unit tangent vector randomly
chosen in the link of the starting point.

In \cite{TB} it is also proved that $Y_t^{\eta}$ (for $\eta\in (0,1]$) is a continuous
Markov process, for each $\eta\in(0,1]$, $Y^\eta$ generates a measure $\mu_\eta$ on the
space $\mathcal{C}(\mathbb{R}^+,K):=\{f:\mathbb{R}^+\rightarrow K, f-continuous\}$ and
$\mu_\eta$ has a subsequence which converges to a measure $W$ on
$\mathcal{C}(\mathbb{R}^+,K)$, called {\em Wiener measure}. This Wiener measure generate a
Brownian motion in the Riemannian polyhedron, such that the transition functions of the
generated Brownian motion are just the projections of the Wiener measure on $K$ (see for
details \cite{TB}).

\begin{prop}\label{brownian}
Let $(B_t)_{t\geq0}$ denote the Brownian motion in the $n-$ dimensional admissible
Riemannian polyhedron $K$, see \cite{TB}. Then $(B_t)_{t\geq0}$ almost surely never hits the
$(n-2)-$skeleton and for every point $p$ of the $((n-1)-skeleton)
\setminus((n-2)-skeleton)$, all the maximal simplexes adjacent to $p$ have the same
probability to be chosen by $B_t^p$ i.e. $(B_t)_{t\geq0}$ has equal branch probabilities.
\end{prop}

\proof The {\em $s-$skeleton} of $K$ is usually denoted by $K^{(s)}$, \cite{EF}.

For any $p\in K^{(n-1)}\setminus K^{(n-2)},$ denote by $U\subset K \setminus K^{(n-2)} $ an
open neighbourhood of $p$  and by $\tau_U:=inf \{t>0/B_t^p \notin U\}$ the fist exit time of
$B_t^p $ from $U$.

For any maximal simplex $S$ of $K$ adjacent to $p$ and for $t$ close to $0$,
$$
P(B_{t\wedge \tau_U}^p\in U\cap S)= \mathop{\mbox{\rm lim}}\limits_{\eta\rightarrow 0}
P(^pY_t^\eta \in U\cap S).
$$

Denote by $\tau_1$ the first stopping-time associated to the process $Y_t^\eta$. Suppose
that $P\{\tau_U\leq\tau_1\}=1$. For any $\eta>0$,
$$
P(^pY_t^\eta \in U\cap S)=E[\chi_{ U\cap S}(^pY_t^\eta)]=
$$
$$=
\int\limits_0^\infty e^{-\left(s+\frac{t}{\eta^2}\right)}\int\limits_{\mathfrak{U}_p(U\cap
S)}\chi_{ U\cap S} \left(^p\Upsilon_{\eta\xi}\left(s+\frac{t}{\eta^2}\right)\right)d
\lambda(\xi)\;ds.
$$

So, for all $\eta > 0,\; $  $P(^pY_t^\eta \in U\cap S)$ depends only on the link
$\mathfrak{U}_p(U\cap S)$ which is independent of the choice of the maximal simplex adjacent
to $p.$ We conclude that $B_p^t$ has equal branch probabilities.

Let us compute $P(B_t\in K^{(n-2)}):$
$$P(B_t\in K^{(n-2)})=\mathop{\mbox{\rm lim}}\limits_{\eta\rightarrow 0}
P(Y_t^\eta \in K^{(n-2)}).$$

For any $\eta >0$,
$$P(Y_t^\eta \in K^{(n-2)})=E[\chi_{K^{(n-2)}}(Y_t^\eta)].$$

Since the process $Y_t^\eta$ is Markov, to compute the above average, we can suppose that
$0< t< \tau_1$ which does not change the result. Then
$$
E[\chi_{K^{(n-2)}}(Y_t^\eta)]
$$
is equal to
$$
\int\limits_0^\infty
e^{-\left(s+\frac{t}{\eta^2}\right)}\int\limits_{\mathfrak{U}_p(K(K^{(n-2)}))}
\chi_{K^{(n-2)}}\left(^p\Upsilon_{\eta\xi}\left(s+\frac{t}{\eta^2}\right)\right)d
\lambda(\xi)\;ds,
$$
where $\mathfrak{U_p}(K(K^{(n-2)}))$ denotes all the vectors in the link $\mathfrak{U}_p(K)$
pointing into the $K^{(n-2)}.$

Since $\lambda(\mathfrak{U}_p(K(K^{(n-2)})))=0,$ we have
$$
\int\limits_{\mathfrak{U}_p(K(K^{(n-2)}))} \chi_{
K^{(n-2)}}\left(^p\Upsilon_{\eta\xi}\left(s+\frac{t}{\eta^2}\right)\right)d \lambda(\xi)=0.
$$

Therefore, we have proved that for any $\eta> 0,$ $P(Y_t^\eta \in K^{(n-2)})=0$ and
consequently $P(B_t\in K^{(n-2)})=0.$

\begin{prop}\label{el}
The Brownian motion $(B_t)_{t\geq0}$ on admissible Riemannian polyhedra has an infinitesimal
generator $L$ defined on a Banach subspace $\mathbf{D}_L$ i.e. for every $f\in
\mathbf{D}_L$,
$$
Lf:=\lim_{t\rightarrow 0} \frac{E[f(B_t)]-f(B_0)}{t}\quad \mbox{ uniformly}.
$$
\end{prop}

\proof Remark that the  Brownian motion $(B_t)_{t\geq0}$ is trajectories continuous
\cite{TB} so it is stochastically continuous. Then by Dynkin's result (see \cite{DY} Theorem
2.3), the existence of $L$ is completely insured.
\endproof

Let $(K,g)$ be  a $n$-dimensional admissible Riemannian polyhedron, endowed with a
continuous simplexwise smooth metric.

\begin{definition}
{\rm Consider $U\subset K$ a domain which meets exactly one $(n-1)$-simplex $S_{n-1}$. Let
$S_n^1,\dots,S_n^k$ denote the $n$-simplexes adjacent to $S_{n-1}$ and $f$ a continuous
function on $U$ which is of class $\mathcal{C}^2$ in each
$\mathop{\widehat{S_n^j}}\limits^o\cap U$ and at least of class $\mathcal{C}^1$ in each
$(\mathop{\widehat{S_n^j}}\limits^o\cup \mathop{\widehat{S_{n-1}}}\limits^o)\cap U$. The
function $f$ is said to be of} zero normal trace condition {\em if and only if:
$$
\sum_{j=1}^k D_jf(x)=0
$$
at almost every point $x$ of $S_{n-1}\cap U$, where $D_jf(x)$ denotes the inner normal
derivative of $f_{|S_n^j\cap U}$ at $x$.}
\end{definition}

\begin{rmk}\label{deltadom}
{\rm The space $\mathbf{D}_L$ contains the space
$$
W_{loc}^{1,2}(K)\bigcap \left\{\begin{array}{l}
                    \mbox{\small function of class } \mathcal{C}^2\;
                        \mbox{\small in the interior of the $n$-simplexes}\\
                          \mbox{\small and the  $(n-1)$-simplexes and of zero normal trace condition} \\
                          \end{array}
                          \right\}
$$
which is denoted by $D_L$.}
\end{rmk}

Let $B=(\Omega, \mathcal{F}^0_t,B_t,\theta_t,P)$ be the $K$-valued Brownian motion
introduced above (see \cite{TB}).

\begin{thm}\label{gen}
Suppose that the exponential
map (Section \ref{exp}) is totally geodesic in each point 
in the topological interior of any $(n-1)$-simplex.
Let $\varphi\in W^{1,2}_{loc}(K)$, i.e. there exists a covering 
$\mathcal{U}$ of $K$ with relatively compact subdomains, such that 
$\varphi$ is of finite energy on each $U\in\mathcal{U}$. Assume that 
$\varphi$ is of class $\mathcal{C}^2$ on the interior of each
$n$-simplex and on the interior of each $(n-1)$-simplex, and of zero 
normal trace condition.
Denote by $\tau_U:=inf \{t>0/B_t \notin U\} $ the first exit time of $B_t$ from $U$. Then we
have:

\begin{equation}\label{eq:generat}
\frac{1}{2}\widetilde{\Delta} \varphi=\frac{\partial}{\partial
t}E[\varphi(B_{t\wedge\tau_U})], \mbox { on $U\backslash((n-2)-skeleton)$. }
\end{equation}
where $E[\varphi(B_{t\wedge\tau_U})]$ is the expectation with respect to $B_{t\wedge\tau_U}$
and
$$
\widetilde{\Delta}=\left\{\begin{array}{ll}
                             \frac{1}{k}\sum\limits_{l=1}^k \Delta_l,\; & \mbox{\small at a point in
                             $(n-1)$-skeleton $\setminus (n-2)$-skeleton}, \\
                             &
                            \Delta_l\; \mbox{\small is the Laplacian in }
                            S_n^l\; \mbox{\small at a boundary point;}\\
                            &\\
                           \mbox {the usual Laplacian,} & \mbox{\small in the interior of each simplex}.\\
                              \end{array}\right.$$

\end{thm}

\proof Let $p\in K.$ There are two cases to investigate:

\medskip

{\bf{\em Case }1}: If $p$ is in the topological interior of some $n$-dimensional simplex,
using  \cite{ME} or \cite{EH}, the relation (\ref{eq:generat}) clearly holds.

\medskip

{\bf{ \em Case} 2}: Let $p$ be in the 
$((n-1)-skeleton)\backslash ((n-2)-skeleton)$.

The idea in this case is to transfer, locally, the computations from 
the polyhedron $K$ to its tangent cone over $p,$  $T_pK$.

Suppose that $p$ is in $\mathop{\widehat{S_{n-1}}}\limits^o$ the topological interior of the
$(n-1)$-simplex $S_{n-1}$. Let $S_n^1,S_n^2,\dots,S_n^k$, $k\geq 2$, denote the
$n$-simplexes adjacent to $S_{n-1}$.
Take $V_0$ a neighbourhood of $0$ in
$T_pK$, and consider the exponential map
$$
E_p : V_0\cap T_pS_n^l\rightarrow  S_n^l,
$$
defined in Section \ref{exp}.

We shrink $V_0$ and $U_p$, if necessary, such that $E_p(V_0) = U_p$. Denote by $\Phi_p : U_p \rightarrow V_0$ the inverse map of $E_p$. 
By hypothesis, the maps $E_p$ and $\Phi_p$, are locally totally
geodesic diffeomorphisms onto their images.

\medskip

Let $(B_t)$ denotes the Brownian motion  in the polyhedron $(K,g)$ (see \cite{TB}) and  by
$X_t= \Phi_p(B_t)$. We remark that $(X_{t\wedge \tau_{U_p}})$ is a Brownian motion in the
flat polyhedron $V_0\bigcap T_pK$ with equal branch probabilities (since $\Phi_p$ is totally
geodesic map and $B_t$ has equal branch probabilities), where $\tau_{U_p}$ is the first exit
time of $B_t$ from $U_p$.

\medskip

Using the fact that $V_0\cap T_pS_n^l$ is also a manifold with boundary, we consider the
boundary normal coordinates $(x_1,\dots,x_{n-1},x_n)$ in the neighbourhood $V_0\cap
T_pS_{n-1}$ of $0$, as follows. Pick up coordinates $(x_1,\dots,x_{n-1})$ on $V_0\cap
T_pS_{n-1}$ and for a point $x$ in the collar neighbourhood of $V_0\cap T_pS_{n-1}$ in
$V_0\cap T_pS_n^l$ define $x_n:=d(x,V_0\cap T_pS_{n-1})$. Remark that this latest choice of
the coordinates chart in the tangent cone over the point $p$ is possible because the metric
$g$ in the polyhedron is continuous.

\medskip

Using the boundary normal coordinates $(x_1,\dots,x_{n-1},x_{n_l})$ in any $(V_0\cap T_p
S_n^l), l=1,\dots,k$, the infinitesimal generator $\Delta^e$ of the Brownian motion $X_t$ on
the flat polyhedron $V_0\cap T_pK$ has the following properties (see \cite{BK}):

\medskip

1) $\Delta^e$ is  defined on the space of:
\begin{itemize}
\item   continuous functions on
the flat polyhedron $V_0$ which are of class $C^2$ in the interior of both $n$-simplexes and
$(n-1)$-simplexes, have continuous second derivatives in the interior of $S_{n-1}$ which are
limits of corresponding directional derivatives from the interior of adjacent faces;

\item  functions with zero normal trace condition for any point in the
($(n-1)$-skeleton)$\backslash$ ($(n-2)$-skeleton), i.e. $\sum \limits_{l=1}^k\frac{\partial
f}{\partial x_{n_l}}=0$;
\end{itemize}

2) $2\Delta^e$ is the usual Laplacian in the interior of each simplex;

3) For $q\in \mathop{\widehat{(V_0\cap T_pS_{n-1})}}\limits^o$, we have
$$
\begin{array}{ccl}
 \Delta^ef(q)&= &\frac{1}{2}\left( \sum\limits_{i=1}^{n-1}
               \frac{\partial^2 f}{\partial {x_i}^2}+\frac{1}{k}\sum\limits_{l=1}^k
               \frac{\partial^2 f}{\partial {x_{n_l}}^2}\right)\\
               &&\\
             &=&\frac{1}{2}  \left[\frac{1}{k}\left(
             \sum\limits_{i=1}^{n-1} \frac{\partial^2 f}{\partial {x_i}^2}+
             \frac{\partial^2 f}{\partial {x_{n_1}}^2}
             \right)+\dots+\right.
          \left.\frac{1}{k}\left( \sum\limits_{i=1}^{n-1}\frac{\partial^2 f}{\partial {x_i}^2}+
          \frac{\partial^2 f}{\partial {x_{n_k}}^2}\right) \right] \\
          &&\\
             &=& \frac{1}{2k}\Delta_1^e f(q)+\dots+\frac{1}{2k}\Delta_k^e f(q)\\
             &&\\
             &=&\frac{1}{2k} \sum\limits_{l=1}^k \Delta_l^e f(q),\\
\end{array}
$$
where
$$
\Delta_l^e f(q)= \sum\limits_{i=1}^{n-1}\frac{\partial^2 f}{\partial {x_i}^2}+
          \frac{\partial^2 f}{\partial {x_{n_l}}^2}
$$
is the usual Laplacian on the manifold $V_0\cap T_pS_n^l$ at a boundary point.

Now, for a point $v_1 \in\mathop{\widehat{(V_0\bigcap T_p S_n^l)}}\limits^o$, the second
order Taylor's development of a function $f$ (as above) in the boundary normal coordinates
has the form:
\begin{equation}\label{Tayleucl}
\begin{array}{rl}
f(v_1)=&f(0)+\partial_{n_l}f(0)x_{n_l}(v_1)
+\frac{1}{2}\partial^2_{n_l}f(0)x^2_{n_l}(v_1)\\
+&\frac{1}{2}\sum\limits_{j=1}^{n-1}\partial_{n_l}\partial_j f(0) x_{n_l}(v_1) x_j(v_1)+
\frac{1}{2}\sum\limits_{i=1}^{n-1}\partial_i\partial_{n_l}f(0) x_i(v_1) x_{n_l}(v_1)\\
+&\sum\limits_{i=1}^{n-1}\partial_i f(0) x_i(v_1)+
\sum\limits_{i,j=1}^{n-1}\partial_i\partial_jf(0) x_i(v_1) x_j(v_1)+o(\varepsilon),\\
\end{array}
\end{equation}
with the notation:
$$
\partial_{n_l}= \frac{\partial}{\partial x_{n_l}},
\partial_i= \frac{\partial}{\partial x_i},
\partial^2_{n_l}= \frac{\partial^2}{\partial {x_{n_l}}^2}.
$$

Using the symmetry of the connection, (\ref{Tayleucl}) reduces to:

\begin{equation}\label{Tayleuclred}
\begin{array}{ccl}
f(v_1)&=&f(0)+\partial_{n_l}f(0)x_{n_l}(v_1)
+\frac{1}{2}\partial^2_{n_l}f(0)x^2_{n_l}(v_1)\\
&&+\sum\limits_{j=1}^{n-1}\partial_j\partial_{n_l} f(0) x_{n_l}(v_1) x_j(v_1)
\\
&&+\sum\limits_{i=1}^{n-1}\partial_i f(0) x_i(v_1)+
\sum\limits_{i,j=1}^{n-1}\partial_i\partial_jf(0) x_i(v_1) x_j(v_1)+o(\varepsilon).
\end{array}
\end{equation}

Equation (\ref{Tayleuclred}) evaluated at $X_{t\wedge\tau_U}$ supposed in
$\mathop{\widehat{(V_0 \bigcap T_pS_n^l)}}\limits^o$ becomes:

\begin{equation}\label{TayleuclBraw}
\begin{array}{rl}
f(X_{t\wedge\tau_U})= & f(0)+\partial_{n_l}f(0)x_{n_l}(X_{t\wedge\tau_U})
+\frac{1}{2}\partial^2_{n_l}f(0)x^2_{n_l}(X_{t\wedge\tau_U})\\
+ & \sum\limits_{j=1}^{n-1}\partial_j\partial_{n_l} f(0) x_{n_l}(X_{t\wedge\tau_U})
x_j(X_{t\wedge\tau_U})\\
+ &\sum\limits_{i=1}^{n-1}\partial_i f(0) x_i(X_{t\wedge\tau_U})+
\sum\limits_{i,j=1}^{n-1}\partial_i\partial_jf(0) x_i(X_{t\wedge\tau_U})
x_j(X_{t\wedge\tau_U})+o(t).
\end{array}
\end{equation}

The process $X_{t\wedge\tau_U}$ can be decomposed into a product $(X^{n-1}_{t\wedge\tau_U},
X_{t\wedge\tau_U}^n)$ of two independent processes such that $X_{t\wedge\tau_U}^{n-1}$ is
$(n-1)-$dimensional Euclidean Brownian motion in the submanifold $V_0\cap T_pS_{n-1}$ and
$X_{t\wedge\tau_U}^n$ is a Brownian motion on a graph $\Gamma$ with $k$ edges
$e_1,\dots,e_k$ of length $\varepsilon$ attached to the point $0$ with equal branch
probabilities at $0$ i.e. $X_{t\wedge\tau_U}^n$ is one dimensional Brownian motion on each
edge of $\Gamma$ with equal branching probabilities, see \cite{BC}.

Then by taking the averages and using \cite{BC}, (\ref{TayleuclBraw}) turns into:

\begin{equation}\label{expecteucl}
\begin{array}{rl}
  E^0[f(X_{t\wedge\tau_U})]= &f(0)+ \sum\limits_{l=1}^k
                \partial_{n_l}f(0)E^0[d(o,X^n_{t\wedge\tau_U}) / \{X^n_{t\wedge\tau_U}\in e_l\}]\\
   + & \sum\limits_{l=1}^k \frac{1}{2}\partial^2_{n_l}f(0)
   E^0[d(o,X^n_{t\wedge\tau_U})^2/ \{X^n_{t\wedge\tau_U}\in e_l\}]\\
   + & \sum\limits_{l=1}^k \sum\limits_{i=1}^{n-1}
   \partial_i\partial_{n_l}f(0)E^0[x_i(X^{n-1}_{t\wedge\tau_U})]E^0[d(o,X^n_{t\wedge\tau_U}) /
   \{X^n_{t\wedge\tau_U}\in e_l\}])\\
   + & \sum\limits_{i=1}^{n-1}\partial_if(0)E^0[x_i(X^{n-1}_{t\wedge\tau_U})] \\
   + & \frac{1}{2}\sum\limits_{i,j=1}^{n-1}\partial_i\partial_jf(0)E^0[x_i(X^{n-1}_{t\wedge\tau_U})
   x_j(X^{n-1}_{t\wedge\tau_U})]+o(t).
\end{array}
\end{equation}

Since the process $X^n_t$ is one dimensional Brownian motion on each edge of $\Gamma$ with
equal branch probabilities (see \cite{BK} or \cite{BC}) then :

$$
E^0[d(o,X^n_{t\wedge\tau_U}) / \{X^n_{t\wedge\tau_U}\in
e_l\}]=\frac{1}{k}\sqrt{\frac{2t}{\pi}}
$$
and
$$
E^0[d(o,X^n_{t\wedge\tau_U})^2 / \{X^n_{t\wedge\tau_U}\in e_l\}]=\frac{t}{k}
$$
for every $l$.

So for a function $f$ with zero normal trace condition, we have:
$$
\sum\limits_{l=1}^k \partial_{n_l}f(0)E^0[d(o,X^n_{t\wedge\tau_U}) /
\{X^n_{t\wedge\tau_U}\in e_l\}] = \frac{1}{k}\sqrt{\frac{2t}{\pi}} \sum\limits_{l=1}^k
\partial_{n_l}f(p)=0,
$$

On the other hand, for the $(n-1)$-Euclidean Brownian motion $X^{n-1}_t$, we have:

$$
E^0[x_i(X^{n-1}_{t\wedge\tau_U})]=0
$$
and
$$
E^0[x_i(X^{n-1}_{t\wedge\tau_U}) x_j(X^{n-1}_{t\wedge\tau_U})]= \delta_{ij} t,
$$
where $\delta_{ij} = 1$ if $i=j$  $\delta_{ij} = 0$ if $i\neq j$.

Then equality (\ref{expecteucl}) reduces to:
\begin{equation}\label{expecteuclfinal1}
\begin{array}{cc}
  E^0[f(X_{t\wedge\tau_U})]=&f(0)+\frac{t}{2k}\sum\limits_{l=1}^k \partial^2 _{n_l} f(0)+
                      \frac{t}{2}\sum\limits_{i,j=1}^{n-1}\delta_{ij} \partial_i\partial_jf(0) +o(t). \\
\end{array}
\end{equation}

Which can be written as:

 \begin{equation}\label{expecteuclfinal2}
\begin{array}{ccl}
  E^0[f(X_{t\wedge\tau_U})]&=&f(0)+
   \frac{t}{2k}\sum\limits_{l=1}^k \Delta_l^e f(0) +o(t) \\
  & =& f(0)+t \Delta^e f(0) +o(t)\\
\end{array}
\end{equation}

We infer from (\ref{expecteuclfinal2}) that:

\begin{equation}\label{DynkinEucl}
\begin{array}{cc}
  f(X_{t\wedge\tau_U})= &f(0)+ \frac{1}{k}\sum\limits_{l=1}^k
  \frac{1}{2}\int\limits_0^{t\wedge\tau_U}
  \Delta_l^e f(X_s) \chi_{\{X_s\in T_pS_n^l\}} ds \\
  & + \textsc{ some local martingale,}\\
\end{array}
\end{equation}
for a function $f$ defined on $V_0$ with zero normal trace condition.

Observe  that the zero normal trace condition is preserved by exponential map.

Now, for a function $\varphi$ of class $W^{1,2}_{loc}$, which is of class $C^2$ in both the
topological interior of the $n$-dimensional faces and the $(n-1)$-dimensional faces of the
polyhedron $K$, the equation (\ref{DynkinEucl}) reads for the function $\varphi\circ E_p$:
\begin{equation}\label{DynkinExp}
\begin{array}{cc}
  \varphi\circ E_p(X_{t\wedge\tau_U})= &\varphi\circ E_p(0)+ \frac{1}{k}
  \sum\limits_{l=1}^k \frac{1}{2}\int\limits_0^{t\wedge\tau_U}
  \Delta_l^e (\varphi\circ E_p) (X_s) \chi_{\{X_s\in T_pS_n^l\}} ds \\
  & + \textsc{ some local martingale}\\
\end{array}
\end{equation}

The process $X_t$ is an Euclidean Brownian motion in each maximal face, so we can write (see
\cite{ME}, Proposition 5.18):
\begin{equation}\label{EmeryExpeucl}
\begin{array}{cc}
  \varphi\circ E_p(X_{t\wedge\tau_U})= &\varphi\circ E_p(0)+ \frac{1}{k}
  \sum\limits_{l=1}^k \frac{1}{2}\int\limits_0^{t\wedge\tau_U}
  \mbox{Hess}_l^e (\varphi\circ E_p) (dX,dX) \chi_{\{X_s\in T_pS_n^l\}}  \\
  & + \textsc{ some local martingale,}\\

  \end{array}
\end{equation}
where $\mbox{Hess}_l^e (\varphi\circ E_p)$ denotes the Euclidean Hessian of the function
$(\varphi\circ E_p)$ on the face $T_pS_n^l$.

Since the map $E_p$ is totally geodesic on each maximal face and on the $(n-1)$-dimensional
face of $V_0\cap T_pK$,  using (4.21) and (4.32) from \cite{ME}, we obtain:
\begin{equation}\label{EmeryExp}
\begin{array}{rl}
\varphi(B_{t\wedge\tau_U})= &\varphi\circ E_p(X_{t\wedge\tau_U}) \\
= &\varphi(p)+ \frac{1}{k}
  \sum\limits_{l=1}^k \frac{1}{2}\int\limits_0^{t\wedge\tau_U}
 (T^*E_p \otimes T^*E_p) \overline{\mbox{Hess}_l}\varphi(dX,dX) \chi_{\{X_s\in T_pS_n^l\}}  \\
  & + \textsc{ some local martingale}\\
= & \varphi(p)+ \frac{1}{k}
  \sum\limits_{l=1}^k \frac{1}{2}\int\limits_0^{t\wedge\tau_U}
  \overline{\mbox{Hess}_l}\varphi(d(E_p\circ X),d(E_p\circ X)) \chi_{\{X_s\in T_pS_n^l\}}  \\
  & + \textsc{ some local martingale,}
\end{array}
\end{equation}
where $\overline{\mbox{Hess}_l}\varphi$ denotes the Hessian of the function $\varphi$ in the
boundary normal coordinates defined by $E_p : V_0\bigcap T_pS_n^l\subset \mathbb{R}^n
\rightarrow U_p\bigcap S_n^l \subset K$ in $S_n^l$ viewed as a manifold with boundary (see
Section \ref{exp}).

Relation (\ref{EmeryExp}) is equivalent to (see \cite{ME} (3.13)):
\begin{equation}\label{EmeryExp1}
\begin{array}{cc}
 \varphi(B_{t\wedge\tau_U}) = &\varphi(p)+ \frac{1}{k}
  \sum\limits_{l=1}^k \frac{1}{2}\int\limits_0^{t\wedge\tau_U}
  \overline{\mbox{Hess}_l} \varphi(dB,dB) \chi_{\{B_s\in U\bigcap S_n^l\}}  \\
  & + \textsc{ some local martingale}\\
\end{array}
\end{equation}

By \cite{ME}, Proposition (5.18), we have:

\begin{equation}\label{EmeryExp2}
\begin{array}{cc}
 \varphi(B_{t\wedge\tau_U}) = &\varphi(p)+ \frac{1}{k}
  \sum\limits_{l=1}^k \frac{1}{2}\int\limits_0^{t\wedge\tau_U}
  \Delta_l \varphi(B_s) \chi_{\{B_s\in U\bigcap S_n^l\}}ds  \\
  & + \textsc{ some local martingale,}\\
\end{array}
\end{equation}
where $\Delta_l$ denotes the Laplace-Beltrami operator computed  by using boundary normal
coordinates in a neighbourhood of $p$ on the manifold with boundary $U_p\bigcap S_n^l$.

This concludes the proof.
\endproof

\begin{rmk}
 {\rm The extra-condition appearing in the
hypothesis, that the exponential map is totally
geodesic, even though restrictive, is realized 
in a certain number of cases. For instance,
it holds if $K$ is one-dimensional (tree).
Hence we obtain a new proof of the main
results in \cite{BK}. Another obvious case is
that of flat metrics.}
\end{rmk}

\begin{rmk}\label{deltaonsnl}
{\rm As we have seen, Remark \ref{deltadom} gives us the space of functions on which is
defined the infinitesimal generator of the Brownian motion on the polyhedron. Moreover from
Theorem \ref{gen}, we conclude that the infinitesimal generator is exactly the
Laplace-Beltrami operator on the interior of each simplex and for a point in the
$(n-1)$-skeleton $\setminus$ $(n-2)$-skeleton it is equal to
$\widetilde{\Delta}=\frac{1}{k}\sum\limits_{l=1}^k \Delta_l$ where $\Delta_l$ is the usual
Laplacian  in $S_n^l$ defined at a boundary point.}
\end{rmk}

\begin{lemma}\label{del}
$L$ is uniquely determined on the space $D_L$.
\end{lemma}

\proof All the functions considered are supposed to be at least of class $\mathcal{C}^2$ in
the interior of each $n$-simplex and each $(n-1)$-simplex.

Consider $f\in W_{loc}^{1,2}(K)$. For every $\psi\in W_c^{1,2}(K)$, by Theorem \ref{gen}, we
have
$$
\frac{1}{2}\int\limits_{K\setminus ((n-2)-skeleton)}\psi \widetilde{\Delta} f d\mu_g
=\int\limits_{K\setminus ((n-2)-skeleton)}\psi Lf d\mu_g.
$$

Consider now an operator $\tilde L$  which is weakly defined on the space $W_{loc}^{1,2}(K)$
by:
$$
 \int\limits_{K}\psi \tilde L f d\mu_g:=
-\frac{1}{2}\int\limits_{K}\langle \nabla \psi, \nabla f\rangle d\mu_g,
$$
for every $\psi\in W_c^{1,2}(K)$.

Indeed, $\tilde L$ is well defined since $W_c^{1,2}(K)$ is a Dirichlet space (see \cite{EF}
page 20, 21 and Proposition 5.1) in the Sobolev  $(1,2)$-norm:
$||u||^2=\int\limits_{K}(u^2+|\nabla u|^2),$ for $u:K\rightarrow \mathbb{R}.$

It is clear that:
$$
\int\limits_{K\setminus ((n-2)-skeleton)}\psi L f d\mu_g= \int\limits_{K\setminus
((n-2)-skeleton)}\psi \tilde L f d\mu_g.
$$

On the other hand, the Brownian motion almost surely never hits the $(n-2)$-skeleton, so
$\tilde L$ is also an infinitesimal generator associate to the transition probability $W_t$
of the Brownian motion.

The transition function $W_t$ associate to the $K$-valued Brownian motion is stochastically
continuous,  so its infinitesimal generator is uniquely determined (see \cite{DY} Lemma 2.2,
Theorem 2.3). We infer that $\tilde L$ is equal to $L$ on the space $D_L$, which concludes
the proof.
\endproof

\begin{thm}\label{procc}
Let $(K,g)$ be an admissible Riemannian polyhedron 
endowed with a simplexwise smooth Riemannian metric and 
$f\in D_L$. As in Theorem \ref{gen}, suppose that the 
exponential map is totally geodesic in each point
in the topological interior of an $(n-1)$-simplex.
Let $(B_t)_{t\geq 0}$ be a $K$-valued Brownian motion and
let $U$ be an open set of $K$ taken as in the hypothesis of Theorem \ref{gen}. Then, for any
$p\in U\cap K\backslash ((n-2)-skeleton)$ the process
$$
C_{t\wedge\tau_U}^{f(p)}=f(B^p_{t\wedge\tau_U})-f(p)-\int\limits_0^{t\wedge \tau_U}
L(f(B_s^p))ds
$$
is a local martingale, where $\tau_U:=inf \{t>0/B_t \notin U\} $ is the first exit time of
$B_t$ from $U$.
\end{thm}

\proof By construction, the Brownian motion $(B_t)_{t\geq 0}$ almost surely never hits the
$(n-2)$-skeleton.

For any $p\in U\backslash((n-2)-skeleton)$, consider the process:
$$
\tilde C_{t\wedge\tau_U}^{f(p)}=\chi_{\{B_t^p\notin ((n-2)-skeleton)\}}
f(B^p_{t\wedge\tau_U})-f(p)
$$
$$-\int\limits_0^{t\wedge \tau_U} \chi_{\{B_s^p\notin
((n-2)-skeleton)\}}L(f(B_s^p))ds ,
$$
where $\chi$ denote the characteristic function.

By Theorem \ref{gen}, $\forall p\in U\backslash ((n-2)-skeleton),$
$$\chi_{\{B_s^p\notin ((n-2)-skeleton)\}}L(f(B_s^p))=
\chi_{\{B_s^p\notin ((n-2)-skeleton)\}} \frac{1}{2}\widetilde{\Delta} f(B_s^p).$$

Taking the expectation, we obtain:
$$
E[\tilde C_{t\wedge\tau_U}^{f(p)}]=E[f(B^p_{t\wedge\tau_U})]-E[f(p)]-
\frac{1}{2}E[\int\limits_0^{t\wedge \tau_U}
 \chi_{\{B_s^p\notin ((n-2)-skeleton)\}} \widetilde{\Delta} f(B_s^p)ds].
$$
Moreover,
$$
\frac{\partial}{\partial t}E[\tilde C_t^{f(p)}]= \frac{\partial}{\partial
t}E[f(B^p_t)]-\frac{1}{2}\widetilde{\Delta} f(p).
$$
Then, using Theorem \ref{gen}, we obtain $\frac{\partial}{\partial t}E[\tilde C_t^{f(p)}]=0,
\forall p\in U\backslash ((n-2)-skeleton).$

Hence, $\tilde C_t^{f(p)}$ is a local martingale. Since $\tilde C_t$ is equal, almost
surely, to the process $C_t$ we conclude that $C_t$ is also a local martingale.
\endproof

\begin{cor}\label{th1}
Let $(K,g)$ be an admissible Riemannian polyhedron endowed with a simplexwise smooth
Riemannian metric, $f\in D_L$ and let $U$ be an open set of $K$  considered as in the
hypothesis of Theorem \ref{gen}. Denote by $\tau_U:=inf \{t>0/B_t \notin U\} $ the first
exit time of $B_t$ from $U$. Then $f$ is harmonic if and only if, for any $p\in U\setminus
((n-2)-skeleton)$, $f(B^p_{t\wedge\tau_U})$ is a local martingale ($(B^p_t)_t$ is a
$K$-valued Brownian motion).
\end{cor}

\proof By Theorem \ref{procc}, the processes:
$$
C_{t\wedge\tau_U}^{f(p)}=f(B^p_{t\wedge\tau_U})-f(p)-\int \limits_0^{t\wedge\tau_U}
L(f(B_s^p))ds
$$
and
$$
\tilde C_{t\wedge\tau_U}^{f(p)}= \chi_{\{B_t^p\notin ((n-2)-skeleton)\}}
f(B^p_{t\wedge\tau_U})-f(p)
$$
$$
- \int\limits_0^{t\wedge\tau_U} \chi_{\{B_s^p\notin ((n-2)-skeleton)\}} L(f(B_s^p))ds
$$
are both local martingales, for every $p\in U\setminus ((n-2)-skeleton)$, where $U$ is taken
as in the hypothesis of the theorem.

Suppose that $f$ is harmonic, then:
$$
\chi_{\{B_t^p\notin ((n-2)-skeleton)\}}f(B^p_{t\wedge \tau_U})= \tilde C^{f(p)}_{t\wedge
\tau_U} +f(p),
$$
for every $p\in U\setminus ((n-2)-skeleton)$.

So the process $\chi_{\{B_t^p\notin ((n-2)-skeleton)\}}f(B^p_{t\wedge \tau_U})$ is a local
martingale, for every $p\in U\setminus ((n-2)-skeleton)$. But this last process is almost
surely equal to $f(B^p_{t\wedge\tau_U})$, so the process $f(B^p_{t\wedge\tau_U})$ is also a
locale martingale.

Conversely, suppose that for every $p\in U\setminus ((n-2)-skeleton)$, $f(B^p_{t\wedge
\tau_U})$ is a local martingale. Then, by classical theory, this implies that $f$ is
harmonic on each $U\setminus ((n-2)-skeleton)$, so is an $E$-minimizer on each $U\setminus
((n-2)-skeleton)$. Then we have for every $\psi\in W_{loc}^{1,2}(K)$, with $\psi=f$ on
$K\setminus U$,
$$
\int\limits_U e(f)=\int\limits_{U\setminus ((n-2)-skeleton)}e(f) \leq\int\limits_{U\setminus
((n-2)-skeleton)}e(\psi)=\int\limits_U e(\psi).
$$

We infer that $f$ is a continuous locally $E$-minimizer map on $K$, which means that $f$ is
harmonic on $U$.
\endproof

\section{Brownian motions, Harmonic maps and morphisms.}

In this section, we extend classical results due to Darling \cite{D} relating Brownian
motions and harmonic maps and morphisms to the case of maps defined on a Riemannian
polyhedron. We prove that harmonic maps are characterized by mapping Brownian motions into
martingales, Theorem \ref{th2}, and harmonic morphisms are exactly the maps which are
Brownian preserving paths, Theorem \ref{th3}.

\begin{thm}\label{th2}
Let $(K,g)$ be an admissible Riemannian polyhedron
as in Theorem \ref{gen} and $(N,h)$ a smooth Riemannian manifold. Let $\varphi:K\rightarrow N$ be
a continuous map such that $\varphi\in W^{1,2}_{loc}(K,N)$, i.e. there exists a covering
$\mathcal{U}$ of $K$ with relatively compact subdomains, such that $\varphi$ is of finite
energy on each $U\in\mathcal{U}$. Assume that, $\varphi$ is  of class $\mathcal{C}^2$ on the
interior of each $n$-simplex and on the interior of each $(n-1)$-simplex. Suppose that for
every function $\psi: N \rightarrow \mathbb{R}$ of class $\mathcal{C}^2$, the function $\psi
\circ \varphi$ is of zero normal trace condition.

Then $\varphi$ is harmonic if and only if for almost all $p\in U$ (with respect to the
volume measure), $\varphi(B^p_{t\wedge\tau_U})$ is a $^N\nabla$-martingale, where
$(B^p_t)_{t\geq 0}$ is a $K$-valued Brownian motion and $\tau_U$ is the first exit time of
$B_t$ from $U$.
\end{thm}

%\begin{rmk}

%Using Whitney's embedding theorem, it is enough to test the
%hypothesis on the map $\varphi$ with a given (any) coordinates
%functions on the Riemannian manifold $N$.
%\end{rmk}

\proof Let $(U_N,V_N,W_N,f)$ be a $^N\nabla$-martingale tester on $N$, such that
$\varphi^{-1}(W_N)\subset U.$

Suppose that $\varphi$ is a harmonic map.

By Theorem \ref{procc}, for all $p\in U\setminus ((n-2)-skeleton)$, the processes:
$$
C_{t\wedge \tau_U}^{f\circ\varphi(p)}:=
f\circ\varphi(B^p_{t\wedge\tau_U})-f\circ\varphi(B^p_0)- \int\limits_0^{t\wedge
\tau_U}\chi_F(s) L(f\circ\varphi)(B^p_s)ds
$$
and
$$
\begin{array}{cc}
  \tilde C_{t\wedge \tau_U}^{f\circ\varphi(p)}:= &
  \chi_{\{B_t^p\notin ((n-2)-skeleton)\}} f\circ\varphi(B^p_{t\wedge\tau_U})
  -f\circ\varphi(B^p_0)- \\
  &\\
   & \int\limits_0^{t\wedge \tau_U}
 \chi_{\{B_s^p\notin ((n-2)-skeleton)\}}\chi_F(s) L(f\circ\varphi)(B^p_s)ds \\
\end{array}
$$
are local martingales, where $F=\mathop\bigcup\limits_{i=1}^{\infty}](\sigma_i,\tau_i]$,
with
$$
\begin{array}{c}
 \sigma_i=inf\{t>\tau_{i-1};\;\varphi( B^p_t)\in U_N\} \\
  \tau_i=inf\{t>\sigma_i;\; \varphi(B^p_t)\not\in \overline{V}_N\}\\
 \sigma_0=0 \\
   \tau_0=0.\\
\end{array}
$$

The map $\varphi$ is supposed to be harmonic and $f$ is a convex function, hence by
Eells-Fuglede's result (\cite{EF}, Theorem 12.1), $f\circ\varphi$ is  a subharmonic function
on $\varphi^{-1}(W_N)$. But in our case the subharmonicity can be translated by:
$$
\widetilde{\Delta}(f\circ\varphi)(p)\geq 0, \forall p\in U\backslash ((n-2)-skeleton)
$$

The process $\chi_{\{B_t^p\notin ((n-2)-skeleton)\}} (f\circ\varphi)(B^p_{t\wedge \tau_U})$
is then the sum of a local martingale and an increasing process, so it is  a local
submartingale $\forall p\in U\backslash ((n-2)-skeleton)$.

Since the process $\chi_{\{B_t^p\notin ((n-2)-skeleton)\}} (f\circ\varphi)(B^p_{t\wedge
\tau_U})$ is equal almost surely  to $(f\circ\varphi)(B^p_{t\wedge \tau_U})$, this last
process is also a local submartingale for every $p\in U\backslash ((n-2)-skeleton)$.

Conversely, suppose that for any $p\in U\setminus ((n-2)-skeleton)$, ($U$ is as in the
hypothesis of the theorem ), $\varphi(B^p_{t\wedge \tau_U})$ is a $^N\nabla$-martingale.

Therefore, for any $^N\nabla$-tester function $f:W_N\rightarrow {\bf R}$,
$(f\circ\varphi)(B^p_{t\wedge \tau_U})$ is a local submartingale.

By Theorem (\ref{procc}), the process:
$$
H^{f\circ\varphi(p)}_{t\wedge\tau_U}:=f\circ\varphi(B^p_{t\wedge \tau_U})-
f\circ\varphi(B^p_0)-\int\limits_0^{t\wedge \tau_U}\chi_F(s) L(f\circ\varphi)(B_s^p)ds,
$$
is a local martingale, for any $p\in U\setminus ((n-2)-skeleton)$, where
$F=\mathop\bigcup\limits_{i=1}^{\infty} (\sigma_i,\tau_i]$, with
$$
\begin{array}{c}
  \sigma_i=inf\{t>\tau_{i-1}; \varphi(B^p_t)\in U_N\} \\
 \tau_i=inf\{t>\sigma_i;\varphi( B^p_t)\not\in \overline{V}_N\}\\
 \sigma_0=0 \\
  \tau_0=0.\\
\end{array}
$$

Since  for any $p\in U\setminus ((n-2)-skeleton)$ we have
$$
\frac{\partial}{\partial t}E[(f\circ\varphi)(B^p_{t\wedge \tau_U})]=\frac{1}{2}
\widetilde{\Delta} (f\circ\varphi)(p).
$$
and $(f\circ\varphi)(B^p_{t\wedge \tau_U})$ is a local submartingale, then
$\widetilde{\Delta} (f\circ\varphi)(p)\geq 0$, for any $p\in U\setminus ((n-2)-skeleton).$

Hence, by Eells-Fuglede (see \cite{EF}, Theorem 12.1), we obtain that $\varphi$ is a
harmonic map on $U\setminus ((n-2)-skeleton).$

Using the same arguments as in the proof of Corollary \ref{th1}, for harmonic functions, we
conclude that $\varphi$ is harmonic on each $U$.
\endproof

\begin{thm}\label{th3}
Notation as in Theorem \ref{th2}. Then $\varphi$ is a 
harmonic morphism if and only if $\varphi$ maps
$K$-valued Brownian motions $(B^p_{t\wedge \tau_U})_{t\geq0}$, for any $p\in U\cap
K\setminus ((n-2)-skeleton)$, to a Brownian motion on $N$, i.e. if $(^NB^p_t)_{t\geq0}$
denote the Brownian motion on the manifold $N$ then there exist a continuous increasing
process $(A_{t\wedge \tau_U})_{t\geq0}$ such that: $^NB^p_t\circ A_{t\wedge \tau_U}=\varphi
\circ B^p_{t\wedge \tau_U}$. $\tau_U$ denote the first exit time of $B_t$ from $U$.
\end{thm}

\begin{rmk}{\rm
We shall suppose the $\dim K\geq \dim N$. Otherwise, $\varphi$ is constant. Indeed, if $\dim
K<\dim N$, by smooth theory (see \cite{BW}, p.46), $\varphi$ is constant on each interior of
maximal simplex (of a chosen fine triangulation of $K$). On the other hand, $\varphi$ is
continuous and $K$ is $(n-1)$-chainable, so $\varphi$ is constant on $K$.}
\end{rmk}

\proofth\ref {th3}: For the proof of the Theorem \ref{th3}, we will adapt and complete the
proof given by Darling (see \cite{D}) in the smooth case.

\medskip

$"\Rightarrow"$\;\; Suppose $\varphi:K\rightarrow N$ is a harmonic morphism. Then by
Eells-Fuglede's result (\cite{EF}, Theorem 13.2), there exists a function $\lambda\in
L^1_{loc}(K)$ called the {\em dilation}, such that:
\begin{equation}\label{dilation}
-\int\limits_K\langle\nabla\psi,\nabla(v\circ\varphi)\rangle=
\int\limits_K\psi\lambda[(\Delta_Nv)\circ\varphi],
\end{equation}
for every $v\in \mathcal{C}^2(N)$ and $\psi\in Lip_c(K)$.

Let $B^p=(B^p_t)_{t\geq0}$ a $K$-valued Brownian motion, for any $p\in K\setminus
((n-2)-skeleton)$. Take $U$ as in the hypothesis of the theorem and suppose $p\in U$.

Define a continuous increasing process $(A_{t\wedge \tau_U})_{t\geq 0}$ by:
\begin{equation}\label{proc}
A_{t\wedge \tau_U}:=\int\limits_0^{t\wedge \tau_U}\lambda(B_s)ds,
\end{equation}
and it's inverse as: $ C_{t\wedge \tau_U}=inf\{s;A_{s\wedge \tau_U}>s\}. $

Denote $\varphi\circ B^p$ by $X^{\varphi(p)}=(X_t^{\varphi(p)})_{t\geq0}$ on $N$. For any
function $f:N\rightarrow \mathbb{R}$ of class $\mathcal{C}^2$ and any $p\in U\cap K\setminus
((n-2)-skeleton)$ we have:
$$
\int\limits_0^{t\wedge \tau_U}\Delta_Nf(X^{\varphi(p)}\circ C_s)ds=
\int\limits_0^{C_{t\wedge \tau_U}}\Delta_Nf(X^{\varphi(p)}_u)dA_u,
$$
where $\Delta_N$ denote the Laplace-Beltrami operator on the manifold $N$.

From (\ref{proc}) we obtain
\begin{equation}\label{inv}
\int\limits_0^{t\wedge \tau_U}\Delta_Nf(X^{\varphi(p)}\circ C_s)ds=
\int\limits_0^{C_{t\wedge \tau_U}}\lambda(B_u^p)\Delta_Nf(X^{\varphi(p)}_u)du.
\end{equation}

Using (\ref{dilation}), the right hand side of the equality (\ref{inv}) is equal to
$$
2\int\limits_0^{C_{t\wedge \tau_U}}L(f\circ\varphi)(B_u^p)du, \; \mu_g \mbox{ a.e. },
$$
for every $p\in U\cap K\setminus ((n-2)-skeleton)$.

On the other hand, by Theorem \ref{procc}, the process
$$
\begin{array}{cc}
  \tilde H_{t\wedge \tau_U}^{f\circ\varphi(p)}:=
    & \chi_{\{B_t^p\notin ((n-2)-skeleton)\}}(f\circ\varphi)(B_{t\wedge \tau_U}^p)-(f\circ\varphi)(p)- \\
   &  \\
   &\int\limits_0^{t\wedge \tau_U} \chi_{\{B_s^p\notin ((n-2)-skeleton)\}} L(f\circ\varphi)(B_s^p)ds\\
\end{array}
$$
is a continuous local martingale for any $p\in U\cap K\setminus ((n-2)-skeleton)$.

Consider now the process
$$\tilde H^{f\circ\varphi(p)}\circ C_{t\wedge \tau_U}\stackrel{denote}
{=}\tilde R_{t\wedge \tau_U}^{f\circ\varphi(p)}.$$ $\tilde R_{t\wedge
\tau_U}^{f\circ\varphi(p)}$ is obviously a continuous local martingale and it is also almost
surely equal (using (\ref{dilation}) and (\ref{inv})) to the process
$$
R_{t\wedge \tau_U}^{f\circ\varphi(p)}:= f(X^{\varphi(p)}\circ C_{t\wedge
\tau_U})-f(\varphi(p))- \frac{1}{2}\int\limits_0^{t\wedge
\tau_U}\Delta_Nf(X^{\varphi(p)}\circ C_s)ds.
$$

So $R_{t\wedge \tau_U}^{f\circ\varphi(p)}$ is also a continuous local martingale for any
$p\in U\cap K\setminus ((n-2)-skeleton)$, which means, by definition, that
$X^{\varphi(p)}\circ C_{s\wedge \tau_U}$ is a Brownian motion on $N$.

\medskip

$``\Leftarrow``$\;\; Conversely, suppose that for any $p\in K\setminus ((n-2)-skeleton)$,
$(\varphi(B_t^p))_{t\geq0}$ is a Brownian motion on $N$ up to a change of time.

Let $V$ be an open set of $N$ such that $\varphi^{-1}(V)\subset U$, where $U$ is taken as in
the hypothesis of the theorem. Let $f:V\rightarrow \mathbb{R}$ be a local harmonic function
on $N$.

Fix $p_0\in U\setminus ((n-2)-skeleton)$ with $\varphi(p_0)\in V$ and $\tau$ denote the
first exit time of $(B^{p_0}_t)_{t\geq 0}$ from $\varphi^{-1}(V)$.

By hypothesis, the process
$$
(f(\varphi(B^{p_0}_{t\wedge\tau})))_{t\geq 0}= (f\circ\varphi(B^{p_0}_{t\wedge\tau}))_{t\geq
0}
$$
is equal to $f(^NB^{\varphi(p)}\circ A_{t\wedge\tau})$.

The latter process is a continuous local martingale (because, by definition the process
$(f(^NB_{t\wedge\tau}^{\varphi(p)}))_{t\geq 0}$ is a continuous local martingale and the
martingale property is stable under change of time).

So we have shown that for every $p_0\in U\setminus ((n-2)-skeleton)$ and for every (local)
harmonic function on $N$, $(f\circ\varphi(B^{p_0}_{t\wedge\tau}))_{t\geq 0}$ is a continuous
local martingale. By Corollary \ref{th1} this means that $f\circ\varphi$ is harmonic.

In other words, we have shown that $\varphi$ pulls-back (local) harmonic function on $N$ to
(local) harmonic function on $K\setminus ((n-2)-skeleton)$. But we have already proved in
the proof of the Corollary \ref{th1} that (local) harmonic function on $K\setminus
((n-2)-skeleton)$ are (local) harmonic on $K$.

We conclude that $\varphi$ is a harmonic morphism.
\endproof

\vskip .5cm

\small \noindent M. A. Aprodu: Department of Mathematics, University of Gala\c ti,
Domneasc\u a Str. 47, RO-6200, Gala\c ti, Romania. e-mail: Monica.Aprodu@ugal.ro
\\\\
T. Bouziane: The Abdus Salam International Centre for Theoretical Physics, strada Costiera
11, 34014 Trieste, Italy. e-mail: taoufik.bouziane@gmail.com

\end{document}